\documentclass[11pt]{amsart}
\input epsf
\usepackage{amsmath}
\usepackage{amssymb}
\usepackage{amscd}
\usepackage{amsthm}
\usepackage{yhmath}
\usepackage{subfigure}

\usepackage[all]{xy}
\usepackage{color}

\numberwithin{equation}{section}

\newtheorem{theorem}[equation]{Theorem}

\newtheorem{proposition}[equation]{Proposition}

\newtheorem{lemma}[equation]{Lemma}

\newtheorem{corollary}[equation]{Corollary}

\theoremstyle{remark}
\newtheorem{remark}[equation]{Remark}

\theoremstyle{definition}

\newtheorem{question}[equation]{Question}

\def\XXint#1#2#3{{\setbox0=\hbox{$#1{#2#3}{\int}$} 
	\vcenter{\hbox{$#2#3$}}\kern-.5\wd0}}

\newcommand{\N}{\mathbb N}

\newcommand{\R}{\mathbb R}

\newcommand{\G}{{\mathbb G}}

\newcommand{\F}{{\mathcal F}}

\newcommand{\M}{{\mathcal M}}

\newcommand{\bilip}{biLipschitz }

\newcommand{\Tan}{\operatorname{Tan}}
\newcommand{\g}{\mathfrak{g}}
\newcommand{\Length}{\operatorname{Length}}

\newcommand{\norm}[1]{\left\Vert#1\right\Vert}

\def\eps{\epsilon}

\begin{document}

\title{Metric spaces with unique tangents}

\author{Enrico Le Donne}

\date{December 9, 2010}

\begin{abstract}
We are interested in studying doubling metric spaces with the property that at some of the points 
the metric tangent is unique. In such a setting, Finsler-Carnot-Carath\'eodory geometries and Carnot groups appear as models for the tangents. 
The results are based on an analogue for metric spaces of 
Preiss's phenomenon: tangents   of tangents are tangents.
\end{abstract}

\maketitle

\section{Introduction}

This paper shows that there is a relation between isometrically homogeneous spaces and  uniqueness of tangents for metric spaces. 
It is a consequence of the work of Gleason, Montgomery-Zippin,  Berestovski\u\i,  Mitchell, and Margulis-Mostow that  a finite-dimensional geodesic metric space 
with transitive  isometry group 
has the property that at every point the tangent metric space is unique. Such a tangent is in fact a  Carnot group equipped with a Carnot-Carath\'eodory distance.
In the following paper we consider doubling-measured   metric spaces with the property that at almost every point the tangent metric space is unique and show that 
almost all tangents have transitive  isometry group. 
Consequently,  if in addition  the metric space  is geodesic, then the tangents are almost surely  Carnot groups equipped with  Carnot-Carath\'eodory distances.

Our  results are founded on the translation in the context of metric spaces of a fact that is well known in Geometric Measure Theory: 
tangent measures of tangent measures are tangent measures.
In fact, let $\mu$  be a doubling measure  in the Euclidean space $\R^n$. Then one can define the tangent measures of $\mu$ at a point $x\in\R^n$ 
by taking weak* limits of measures translated by $x$ and dilated by larger and larger factors. Namely, setting $T_{x,\rho}(y):=\rho(y-x)$, one defines
$$\nu\in \Tan(\mu,x) \iff \nu=\lim_{i\to\infty} c_i(T_{x,\rho_i})_\# \mu , \text{ for some } \rho_i\to\infty \text{ and }  c_i\to 0.$$
In \cite{Preiss}, David Preiss showed the useful fact that, for $\mu$-almost every $x$, 
any tangent measure of a tangent measure of $\mu$ at $x$ is itself a tangent measure of $\mu$ at $x$. 

We shall consider tangents in the class of doubling metric spaces.
Let $(X,d)$ be a metric space whose distance is doubling. Mikhail
Gromov showed that one can consider the tangent spaces at a point $x\in X$ as
 the limits of sequences of   pointed metric spaces $(X, \rho_j d, x)$, with $\rho_j\to\infty$, as 
$j\to \infty$.

We shall consider doubling-measured metric spaces, i.e., metric spaces endowed with a doubling measure. As a consequence, the distance itself is doubling.
We show that the analogue of Preiss's phenomenon holds:
\begin{theorem}\label{tantan}
Let $(X,\mu, d)$ be a doubling-measured metric space. Then the following two properties hold.
\begin{enumerate}
\item For $\mu$-almost every $x\in X$, for all $(Y,y)\in \Tan(X,x)$, and for all $y'\in Y$ we have
$(Y,y')\in\Tan(X,x).$
\item For $\mu$-almost every $x\in X$, for all $(Y,y)\in \Tan(X,x)$, and for all $y'\in Y$ we have
$\Tan(Y,y')\subseteq\Tan(X,x).$
\end{enumerate}
\end{theorem}

Notice that as pointed metric spaces $(Y,y)$ and $(Y,y')$ might be different. A  limit space is defined up to isometry. Hence,  $(Y,y)$ and $(Y,y')$ are equal when there exists an  isometry $f:Y\to Y$ with the property that $f(y)=y'$.
Therefore, if it is  the case that there is only one tangent metric space at a point $x$ where the conclusion of the part (1) of  Theorem \ref{tantan} holds,
then such a metric space 
$(Y,y)$ in $ \Tan(X,x)$ has the property that $(Y,y)$ is isometric to $(Y,y')$, for all $y'\in Y$.
In other words, the isometry group of $Y$ acts on $Y$ transitively.
In conclusion, uniqueness of tangent spaces leads to isometric homogeneity of such tangents.

In the next theorem we completely characterize the tangents that can appear, if in addition  the metric space is geodesic. For more general results see Section 2. 

\begin{theorem} \label{main thm}
Let $(X,\mu, d)$ be a doubling-measured geodesic metric space. Assume that, for $\mu$-almost every $x\in X$, the set $\Tan(X,x)$ contains only one element.
 Then,  for $\mu$-almost every $x\in X$, the element in $\Tan(X,x)$ is a Carnot group $\G$ endowed
 with a sub-Finsler left-invariant metric with the first layer of the Lie algebra of $\G$ as horizontal distribution.
\end{theorem}

Recall that a {\em Carnot group} $\G$ of step $s\geq 1$ is a connected, simply-connected Lie
group whose Lie algebra $\g$ admits a step $s$ stratification: this
means that we can write
$$\g= V_1\oplus \cdots\oplus V_s, $$
with $[V_j, V_1] = V_{j+1}$, for $1\leq j\leq s-1$, and $V_s\neq \{0\}$. The subspace $V_1$ is called the {\em first layer of the Lie algebra} $\g$.
A {\em sub-Finsler left-invariant metric with $V_1$ as horizontal distribution} is defined as follows.
One fixes a norm  $ \norm{\cdot} $ on $V_1$. The space $V_1$ defines a left-invariant sub-bundle of the tangent bundle of $\G$. 
The norm  $ \norm{\cdot} $ extends left-invariantly as well.
The triple $(\G, \Delta,\norm{\cdot} )$ is a left-invariant  sub-Finsler structure for which the 
 \textit{Finsler-Carnot-Carath\'eodory} or \textit{sub-Finsler} distance  $d_{CC}$ is defined as,  for any $x,y\in \G$,
\begin{equation}\label{dist_CC}
d_{CC}(x,y):=\inf\{\Length_ {\norm{\cdot}}(\gamma)\;|\;\gamma\in C^\infty([0,1];\G), \gamma(0)= x , \gamma(1)= y, \dot \gamma\in V_1\}.
\end{equation}

To conclude the introduction, we would like to mention a similar result  of Pertti Mattila, which as well was obtained by proving a  Preiss's phenomenon for measures 
on locally compact groups with metric dilations. 
Namely, in \cite{Mattila-unique} it is shown that if a measure on such a group  has unique tangents, 
then its tangents are almost surely Haar measures with respect to some closed dilation-invariant subgroup.
 
\subsection{Other consequences}

Given a metric space $(X,d)$, we denote the dilated space by a factor $\rho>0$ by
$$\rho X:=\left(X,\rho d\right).$$
Fixed a point $x\in X$, we denote by $\Tan(X,x)$ the space of all Gromov-Hausdorff limits of the pointed metric spaces
$$(\rho_i X, x) ,\qquad \text{ with } \rho_i\to\infty,\qquad \text{ as } i\to \infty.$$
In the next section we will recall   the definition of  Gromov-Hausdorff convergence.
However, consider that the elements in $\Tan(X,x)$ are defined up to isometric equivalence.

A measure $\mu$ on a metric space $(X,d)$ is said to be {\em doubling}, if
there exists a constant $C$ such that, for all $x\in X$ and $r>0$,
$$0\neq\mu(B(x,2r))< C\mu (B(x,r)).$$
Notice that if $\mu$ is a doubling measure, then $d$ is a  {\em doubling distance}, i.e.,   there
is a constant $N$ such that any ball can be
covered by  $N$ balls of half the radius.
Mikhail Gromov showed that, whenever $(X,d)$ is a doubling metric space, then, for any $x\in X$, the set $\Tan(X,x)$ is non-empty.


The following theorem  is a more detailed version of Theorem \ref{main thm}.
For its proof we will use the work of Gleason-Montgomery-Zippin \cite{mz} and the applications  by Berestovski\u\i\; \cite{b, b1, b2}.
\begin{theorem}\label{thm-uniq-tan}
Let $(X,\mu, d)$ be a doubling-measured metric space. 
Let $\Omega\subseteq X$ be the subset of elements 
 $x\in X$ such that the set $\Tan(X,x)$ contains only one element. 
 Then, for $\mu$-almost every $x\in \Omega$, the element in $\Tan(X,x)$ 
is an isometrically homogeneous space of the following form. There is a Lie group $G$ and a compact subgroup $H<G$, such that the tangent space at the point $x$ is isometric to  the  manifold 
 $G/H$ equipped with some $G$-invariant distance function.

If, moreover, the distance $d$ is geodesic, then, for $\mu$-almost every $x\in \Omega$, the element in $\Tan(X,x)$ is a Carnot group $\G$ endowed with a sub-Finsler left-invariant metric with the first layer of the Lie algebra of $\G$ as horizontal distribution.
\end{theorem}

The following is an application of the previous techniques (i.e., Theorem \ref{tantan}) to the theory of \bilip homogeneous spaces.
A metric space $(X,d)$ is said to be {\em locally  biLipschitz  homogeneous} 
if, for every two points
$x_1,x_2 \in X$, there are neighborhoods $U_1$ and $U_2$  of $x_1$ and $x_2$ respectively
and a \bilip homeomorphism $f:U_1\to U_2$, such that $f(x_1)=x_2$.
 We call {\em finite-dimensional isometrically-homogeneous space} a manifold of the form $G/H$, with $G$ a Lie group and $H$ a compact subgroup $H<G$,
 endowed with some $G$-invariant distance function.
\begin{theorem}\label{thm-bilip-homog}
Let $(X,\mu, d)$ be a doubling-measured metric space. Assume that $(X,d)$ is locally biLipschitz homogeneous. Assume also that there are a point $x_0\in X$ and, for some $K>1$, a family of $K$-biLipschitz maps $\F$ such that $\F$ is a group, i.e., $\F$ is closed under composition, and, for all pair of tangents
$(Z_1,z_1), (Z_2,z_2)\in\Tan(X,x_0)$, there is a map $\psi\in\F$ with
$\psi:(Z_1,z_1)\to (Z_2,z_2)$.

Then, there is a finite dimensional isometrically homogeneous space $G/H$, such that, for all $x\in X$, each element in $\Tan(X,x)$ 
is \bilip equivalent to $G/H$.
\end{theorem}

Since the isometries form a group, we immediately have the following result.
\begin{corollary}[of Theorem \ref{thm-bilip-homog}]
Let $(X,\mu, d)$ be a doubling-measured metric space. Assume that $(X,d)$ is locally biLipschitz homogeneous and that, for some point $x_0\in X$ the collection $\Tan(X,x_0)$ contains only one metric space, up to isometric equivalence.

Then,   there is a finite dimensional isometrically homogeneous space $G/H$, such that,  for all $x\in X$, each element in $\Tan(X,x)$ 
is \bilip equivalent to  $G/H$.
\end{corollary}

Next, one should wonder which are the consequences of having tangents equal to Carnot groups. The answer is definitely not easy at least because there are uncountably many Carnot groups.
The case when the tangents are Euclidean is  relatively easier. 
Indeed, G. David and T. Toro considered such a case in their study of Reifenberg flat metric spaces, c.f.~\cite{David-Toro}. 
In the particular case when the metric space is nicely embedded in a Hilbert space, then we observe that the standard `cone criterion' gives the following easy fact.
\begin{corollary}
Let $X$ be a locally compact subset of a separable Hilbert space $H$. Let $d$ be the distance function on $H$ restricted to $X$. Let  $\mu$ be a doubling measure for $(X,d)$.
Assume that, at $\mu$-almost every point $x\in X$, the dilated spaces
$$\rho(X-x)$$
converge in the Hausdorff sense, as $\rho\to\infty$.

Then $X$ is contained in the union of countably many Lipschitz graphs, up to a set of $\mu$-measure $0$.
\end{corollary}

\begin{question} Which are the metric spaces that can be embedded in a separable Hilbert space having the property that, almost everywhere, the Gromov tangents can be calculated as Hausdorff tangents?
\end{question}
\setcounter{tocdepth}{2}
\tableofcontents
\section{Tangents as limit of pointed spaces}
A pointed metric space $(X_\infty, d_\infty, x_\infty)$ is a {\em tangent} of a metric space $(X,d)$ at the point $x\in X$, if there are 
Hausdorff approximations
 $$\left\{\phi_i : (X_\infty, d_\infty, x_\infty) \to (X,d_i,x)\right\}_{i\in\N},$$ with $d_i=\frac{1}{\lambda_i}d$ for some $\lambda_i\to0$ as $i\to\infty$.
 More explicitly, one has that, for all $R\geq0$ and all $\delta>0$,
$$\limsup_{i\to\infty}\left\{|d_i(\phi_i(y), \phi_i(z)) - d_\infty(y, z)| : y, z \in B(x_\infty,R) \subset X_\infty\right\}= 0 $$
and
 $$\limsup_{i\to\infty}\left\{
d_i(y, \phi_i(B(x_\infty,R + \delta))) : y \in B_{d_i}(x,R) \subseteq (X,d_i)
\right\}= 0 .$$

The first condition says that 
$$\dfrac{1}{\lambda_i}d(\phi_i(y), \phi_i(z)) \to d_\infty(y, z),$$
uniformly in $y$ and $z$ on bounded sets.
The second condition can be written as
 $$\limsup_{i\to\infty}\left\{
\dfrac{1}{\lambda_i}d(y, \phi_i(B(x_\infty,R + \delta))) : y \in B(x,\lambda_iR)
\right\}= 0 .$$
Roughly speaking, this means that the sequence of (smaller and smaller) balls $\phi_i(B(x_\infty,R + \delta))$ covers $B(x,\lambda_iR)$ better and better, in fact with a sub-linear gap.

\section{Proofs of the results}

Given a Radon measure $\mu$ on a space $X$, one can consider  the outer measure $\mu^*$, defined for any $A\subseteq X$
by
$$
\mu^*(A):=\inf\left\{\mu(B)\,:\,\text{$B$ Borel, $B\supseteq
A$}\right\}.
$$

Even if any geometric intuition says that the following fact is  obvious, it is the key point  implying Theorem \ref{tantan}.
\begin{proposition}\label{ptdensity}
Let $(X,\mu, d)$ be a doubling-measured metric space. 
Let $A\subseteq X$ be any set and let $a\in A$ be a point of density for $A$, i.e., 
$$\lim_{r\downarrow 0}\frac{\mu^*(A\cap B_r(a))}{\mu(B_r(a))}=1.$$
Then $\Tan(A,d,a)=\Tan(X,d,a).$
\end{proposition}

\begin{corollary}[of the proof of Proposition \ref{ptdensity}]\label{cor:image}
If $(Y,y)\in \Tan(X,x)$ and $x$ is  a point of density for a set $A$, then there are Hausdorff approximations
$\phi_i:Y\to X$ such that ${\rm Im}(\phi_i)\subseteq A$.
\end{corollary}

\begin{lemma}\label{lemmaptdensity}
Let $(X,\mu, d)$ be a doubling-measured metric space. 
Let $A\subseteq X$ be any set and let $a\in A$ be a point of density for $A$. Then, for all $\delta>0$ and all $R\geq 0$, we have
$$\lim_{\lambda\to0}\sup\left\{\dfrac{1}{\lambda} d(p,A\cap B(a,(R+\delta)\lambda))\;:\;p\in B(a, \lambda R)\right\}=0.$$
\end{lemma}

\proof Fix $\eps>0$. 
Let $C$ and $Q$ be the constants of the doubling property for $\mu$, i.e., for all $R>r>0$,
$$\dfrac{\mu(B_R)}{\mu(B_r)}< C\left(\dfrac{R}{r}\right)^{-Q}.$$
Take $\lambda$ small enough such that
$$ \frac{\mu^*( B(a, (R+\delta)\lambda ) \setminus A)}{\mu(B(a,(R+\delta)\lambda))}<
\alpha:=\dfrac{1}{2C}\left(\dfrac{\eps/2}{2R+\delta}\right)^Q.$$
We shall prove that, for such a $\lambda$,
the supremum is smaller than $\eps$. We can assume $\eps/2<\delta$.
Assume, by the way of contradiction, that there is some $p \in B(a, \lambda R)$ such that 
$$d(p,A\cap B(a,(R+\delta)\lambda))\geq \dfrac{\eps}{2}\lambda.$$
Note that, by triangle inequality, we have
$$B(p, \dfrac{\eps}{2}\lambda)\subset B(a,(R+\delta)\lambda)$$
and thus
$$A\cap B(p, \dfrac{\eps}{2}\lambda)=\emptyset.$$
Therefore,
\begin{eqnarray*}
\mu( B(p, \dfrac{\eps}{2}\lambda) ) &\leq & \mu^*( B(a, (R+\delta)\lambda ) \setminus A)\\
&\leq& \alpha \mu(B(a,(R+\delta)\lambda))\\
&\leq& \alpha \mu(B(p,(2R+\delta)\lambda))\\
&\leq& \alpha  C\left(\dfrac{\eps/2}{2R+\delta}\right)^{-Q} \mu (B(p, \dfrac{\eps}{2}\lambda))\\
&\leq& \dfrac{1}{2} \mu (B(p, \dfrac{\eps}{2}\lambda)).
\end{eqnarray*}
This last calculation implies that $1\leq 1/2$, which is a contradiction.\qed

\subsection{Proof of Proposition \ref{ptdensity}}
\proof[Proof of {$\Tan(X,a)\subseteq \Tan(A,a)$}]

Let $(X_\infty, d_\infty, x_\infty)\in \Tan(X,a)$. So there are
Hausdorff approximations 
$$\left\{\phi_i : (X_\infty, d_\infty, x_\infty) \to (X,\dfrac{1}{\lambda_i}d,a)\right\}_{i\in\N},$$ with rescale factors $\lambda_i\to0$.
For all $p\in X_\infty$, define $\phi'_i(p)$ as a closest point in $A$ to $\phi_i(p)$.
Notice that $A$ might be considered closed, since $A$, the completion $\bar{A}$ of $A$, and the closure ${ \mathcal{C}}(A)$ of $A$ in $X$ have Gromov-Hausdorff distance $0$, therefore they have the same tangents:
$$\Tan(A,a)=\Tan(\bar{A},a)=\Tan({\mathcal{C}}(A),a).$$
So, we constructed maps
$$\phi'_i : X_\infty \to A.$$
We claim the following:
\begin{description}
\item[Claim 1] $\dfrac{1}{\lambda_i} d(\phi'_i(\cdot),\phi_i(\cdot))\to 0$ uniformly  on bounded sets,
\item[Claim 2] The maps $$\left\{\phi'_i : (X_\infty, d_\infty, x_\infty) \to (A,\dfrac{1}{\lambda_i}d,a)\right\}_{i\in\N}$$ 
are Hausdorff approximations, and so $(X_\infty, d_\infty, x_\infty)\in \Tan(A,a)$.
\end{description}

\proof[Proof of Claim 1]
Fix $R>0$. 
Observe that, clearly, for any $\delta>0$,
$$d(\cdot, A)\leq d(\cdot, A\cap B(a,(R+\delta)\lambda_i)).$$
Therefore, Lemma \ref{lemmaptdensity} gives
$$\lim_{i\to\infty}\sup\left\{\dfrac{1}{\lambda_i} d(p,A)\;:\;p\in B(a, \lambda_i R)\right\}=0.$$
Fix some $\eta>0$. For $i$ big enough, we have that, for all $q\in B(x_\infty, R-\eta)$,
$$\dfrac{1}{\lambda_i} d(\phi_i(q),a)\leq R.$$
By definition of $\phi'_i$, we have
$$d(\phi_i(q),\phi'_i(q))=d(\phi_i(q),A).$$
Thus,
$$\lim_{i\to\infty}\sup\left\{\dfrac{1}{\lambda_i} d(\phi_i(q),\phi'_i(q))\;:\;q\in B(x_\infty, R-\eta)\right\}=0.$$
\qed
\proof[Proof of Claim 2]
First,
\begin{eqnarray*}
&&
\sup \left\{|d_i(\phi'_i(y), \phi'_i(z)) - d_\infty(y, z)| : y, z \in B(x_\infty,R) \subset X_\infty\right\}\\
&\leq&
\sup\left\{|d_i(\phi_i(y), \phi_i(z)) - d_\infty(y, z)|+
|d_i(\phi'_i(y), \phi_i(y)) |+
|d_i(\phi_i(z), \phi'_i(z)) | : y, z \in B(x_\infty,R)  \right\}\\
&\leq&  
\sup\left\{|d_i(\phi_i(y), \phi_i(z)) - d_\infty(y, z)|
: y, z \in B(x_\infty,R)  \right\}+\\&&\qquad\qquad+
\sup\left\{
|d_i(\phi'_i(y), \phi_i(y)) |
: y, z \in B(x_\infty,R)  \right\}+\\&&\qquad\qquad\quad+
\sup\left\{
|d_i(\phi_i(z), \phi'_i(z)) | : y, z \in B(x_\infty,R)  \right\} \\
&&\to 0 +0+0=0.
\end{eqnarray*}

Second,
 \begin{eqnarray*}
&&\sup\left\{
d_i(y, \phi'_i(B(x_\infty,R + \delta))) : y \in B(a,\lambda_i R)\cap A  
\right\}\\
&\leq&\sup\left\{
d_i(y, \phi'_i(B(x_\infty,R + \delta))) : y \in B(a,\lambda_i R)  
\right\}\\
&\leq&
\sup\left\{
d_i(y, \phi_i(B(x_\infty,R + \delta))) : y \in B(a,\lambda_i R)  
\right\}+\\&&\qquad\qquad\quad+
\sup\left\{
d_i(\phi_i(z), \phi'_i(B(x_\infty,R + \delta))) : z \in B(a,\lambda_i R)  
\right\}\\
&\leq&
\sup\left\{
d_i(y, \phi_i(B(x_\infty,R + \delta))) : y \in B(a,\lambda_i R)  
\right\}+\\&&\qquad\qquad\quad+
\sup\left\{
d_i(\phi_i(z), \phi'_i(z)) : z \in B(a,\lambda_i R)  
\right\}\\
&\to& 0 .
\end{eqnarray*}
\qed

\proof[Proof of {$\Tan(A,a)\subseteq \Tan(X,a)$}]
Vice versa, an element $(X_\infty, d_\infty, x_\infty)\in \Tan(A,a)$ gives
Hausdorff approximations 
$$\left\{\phi_i : (X_\infty, d_\infty, x_\infty) \to (A,\frac{1}{\lambda_i}d,a)\right\}_{i\in\N},$$ with rescale factors $\lambda_i\to0$.

We claim that the following  maps are Hausdorff approximations:
$$\left\{\phi'_i : (X_\infty, d_\infty, x_\infty) \to (X,\frac{1}{\lambda_i}d,a)\right\}_{i\in\N},$$ 
defined as
$$
\phi_i':=\iota\circ \phi_i,$$
where $\iota:A\to X$ is the inclusion.
The first requirement to check is that
$$d_i(\phi'_i(y), \phi'_i(z))=d_i(\phi_i(y), \phi_i(z)) \to d_\infty(y, z),$$
uniformly in $y$ and $z$ on bounded sets, which is clearly true.
The second condition is consequence of Lemma \ref{lemmaptdensity}:
\begin{eqnarray*}
&&\limsup_{i\to\infty}\left\{
d_i(y, \phi'_i(B(x_\infty,R + \delta))) : y \in B(a,\lambda_iR)
\right\}\\
&\leq&\limsup_{i\to\infty}\left\{
d_i(y, A\cap B(a,R + \delta/2))) : y \in B(a,\lambda_iR)
\right\}+\\
&&\qquad+\limsup_{i\to\infty}\left\{
d_i(y, \phi_i(B(x_\infty,R + \delta))) : y \in B(a,\lambda_i(R+\delta/2))\cap A
\right\}\\
&=& 0 .\end{eqnarray*}
\qed

\begin{remark}
Both the doubling property and the density point property are necessary for both the containments.
\end{remark}

\proof[Proof of Corollary \ref{cor:image}]
By Proposition \ref{ptdensity}, the space $(Y,y)$ is in $\Tan(A,x)$. So there are Hausdorff approximations
$$\psi_i :Y\to A.$$
Consider then the maps $\phi_i=\iota \circ \psi_i$, where $\iota:A\to X$ is the inclusion. The calculations at the end of the proof of Proposition \ref{ptdensity} show that such $\phi_i$'s
are Hausdorff approximations for $X$ with image in $A$.
\qed
\subsection{Some facts on the space of metric spaces}
Let $\M$ be the set of all the separable, locally uniformly bounded and pointed metric spaces.
We consider the pointed Gromov-Hausdorff convergence on the set $\M$.

A first fact to recall is that such a topology is metrizable: There exists a distance function $\underline{d}$ on $\M$ such that
$$
(X_\infty, d_\infty, x_\infty)\in \Tan(X,d,x)
\iff
\lim_{\lambda_i\to0} \underline{d}\left( (X, \frac{1}{\lambda_i}d, x),(X_\infty, d_\infty, x_\infty) \right)
=0.$$

A second fact to recall is that the space $(\M,\underline{d})$ is separable. In particular, for each $k\in \N$ there exists a cover 
\begin{equation}\M=\bigcup_{l\in\N}B_l,   \label{cover}
\end{equation}
such that, if $(Y,y)$ and $(Y',y')$ are both in $B_l$,
then 
$$	\underline{d}\left( (Y,y),(Y',y') \right)
< \dfrac{1}{2k}.$$
\subsection{Proof of Theorem \ref{tantan}}
\proof[Proof of Theorem \ref{tantan}]
We need to show that 
$$	\mu\left( \left\{ x\in X : \forall (Y,y)\in\Tan(X,x), \forall y'\in Y: (Y,y')\in\Tan(X,x)\right\}^c  \right)=0.$$
In other words,
$$	\mu\left( \left\{ x\in X : \exists (Y,y)\in\Tan(X,x), \exists y'\in Y: (Y,y')\notin\Tan(X,x)\right\} \right)=0.$$
Using the distance $\underline{d}$ on the collection of metric spaces, we just need to show that, for all $k,m\in\N$, we have
$$	\mu\left( \left\{ x\in X : \exists (Y,y)\in\Tan(X,x), \exists y'\in Y: \underline{d}\left((Y,y'),(\frac{1}{t} X,x)\right)>\frac{1}{k}, \forall t\in (0,\frac{1}{m})\right\} \right)=0.$$
Using the cover \eqref{cover} coming from the separability, we need to show that, for all $k,l,m\in\N$, each set
$$\left\{ x\in X : \exists (Y,y)\in\Tan(X,x), \exists y'\in Y: \right.\qquad\qquad\qquad\qquad\qquad\qquad\qquad$$ 
$$\qquad\qquad\qquad\qquad\qquad\qquad\left.(Y,y')\in B_l\text{ and }\underline{d}\left((Y,y'),(\frac{1}{t} X,x)\right)>\frac{1}{k}, \forall t\in (0,\frac{1}{m})\right\}$$
is $\mu$-negligible.
Assume that one of these sets above is not $\mu$-negligible and call it $A$; so $k,l$ and $m$ are now fixed and $\mu^*(A)>0$. Here we use the outer measure $\mu^*$, since we don't want, and don't need, to show measurability of such a set.

Let $a$ be a point of density of $A$ for $\mu^*$.
Since $a\in A$, there exist $(Y,y)\in\Tan(X,a)$ and $ y'\in Y$ such that $ (Y,y')\in B_l $ and $\underline{d}\left((Y,y'),(\frac{1}{t} X,a)\right)>\frac{1}{k}$, for all $t\in (0,\frac{1}{m}).$

Since $(Y,y)\in\Tan(X,a)$, there is a sequence $\lambda_i\to0$ such that 
$$(\frac{1}{\lambda_i}X,a)\stackrel{GH}{\to}(Y,y).$$
Let $\phi_i:Y\to X$  the Hausdorff approximations with ${\rm Im}(\phi_i)\subseteq A$, given by Corollary \ref{cor:image}.
Let $a_i=\phi_i(y')\in A$. 
Then 
$$(\frac{1}{\lambda_i}X,a_i)\stackrel{GH}{\to}(Y,y').$$
Now take $i$ big enough so that
$$\underline{d}\left(  (\frac{1}{\lambda_i}X,a_i), (Y,y')  \right)<\frac{1}{2k}.$$
Since $a_i\in A$, there are spaces $(Y_i,y_i)\in\Tan(X,a_i)$ and $ y_i'\in Y_i$ such that 
$ (Y_i,y_i')\in B_l $ and $\underline{d}\left((Y_i,y_i'),(\frac{1}{t} X,a_i)\right)>\frac{1}{k}$, for all $t\in (0,\frac{1}{m}).$
So we arrive at a contradiction:
\begin{eqnarray*}
\frac{1}{k}&<&   \underline{d}\left((Y_i,y_i'),(\frac{1}{t} X,a_i)  \right)        \\
&\leq& \underline{d}\left(  (Y_i,y'_i), (Y,y')  \right)+\underline{d}\left(  (\frac{1}{\lambda_i}X,a_i), (Y,y')  \right)\\
&\leq& {\rm Diam}_{\underline{d}}(B_l)  +\underline{d}\left(  (\frac{1}{\lambda_i}X,a_i), (Y,y')  \right)
< \frac{1}{2k}+\frac{1}{2k}.
\end{eqnarray*}
\qed

\subsection{Proof of Theorem \ref{thm-uniq-tan}}
\proof[Proof of Theorem \ref{thm-uniq-tan}]
Using Theorem \ref{tantan}, up to removing a $\mu$-negligible set, we may assume that 
 for all $x\in \Omega$ we have that 
$$\Tan(X,x)=\{(Y,y)\}$$
and for all $y'\in Y$ we also have $(Y,y')\in \Tan(X,x)$.  Thus,
$(Y,y')$ is isometric to $(Y,y'')$ for all $y'$ and $y''\in Y$. 
In other words, the metric space $Y$ is isometrically homogeneous.

Since $Y$ is a tangent of a doubling space, $Y$ is doubling as well.
In particular, the Hausdorff dimension and thus the topological dimension of $Y$ are finite.

By Montgomery-Zippin's work \cite{mz}, the group of isometries $G$ of $Y$ is a Lie group, and $Y$ is homeomorphic to a quotient $G/H$, where $H$ is the stabilizer of a point. Thus there is a $G$-invariant distance on $G/H$ for which $Y$ is isomorphic $G/H$. 

If, moreover, $X$ is geodesic, then $Y$ and $G/H$ are geodesic as well.
 By Berestovskii's Theorem \cite{b}, the $G$-invariant distance function on $G/H$ is a $G$-invariant sub-Finsler metric $d_{SF}$, i.e., there is a $G$-invariant sub-bundle $\Delta$ on the manifold $G/H$ and a $G$-invariant norm on $\Delta$, such that $d_{SF}$ is the Finsler-Carnot-Carath\'eodory distance associated.

We show now that in fact the space $G/H$ is a Carnot group. Indeed, notice  that, since $Y$ is a tangent then any of its dilated spaces is still a tangent. By uniqueness of tangents, such dilations are isomorphic to the same space. In other words, $Y$ is a cone. Consequently, 
$$\Tan(Y,y)\subseteq\Tan(X,x).$$
First, by Mitchell's Theorem \cite{Mitchell}, the tangent to $G/H$ is a Carnot group $\G$. Second, by uniqueness of tangents, we have that $Y=\G$.
\qed
\subsection{Proof of Theorem \ref{thm-bilip-homog}}
\proof[Proof of Theorem \ref{thm-bilip-homog}]
Let $\Omega\subset X$ be a full-measure set for which the conclusion of Theorem \ref{tantan} holds.
Let $x\in \Omega$. Fix any $(Y,y)\in \Tan(X,x)$.

Since $(X,d)$ is \bilip homogeneous, there exists an $L$-\bilip map 
$$f:(U_x,x)\to (U_{x_0},x_0),$$
where $U_x$ and $U_{x_0}$ are neighborhoods of $x$ and $x_0$ respectively.
Since $x\in\Omega$, we have, for all $y'\in Y$, that $(Y,y')\in\Tan(X,x)$.

Let $\lambda_i\to0$ the rescaling factors giving the tangent $(Y,y')$. Consider now the same dilations for the set $X$ but now pointed at $x_0$:
$$(X,\dfrac{1}{\lambda_i}d, x_0).$$ Up to considering a subsequence, since such dilations are uniformly doubling, the sequence converges to a metric space
\begin{equation}\label{zyprime}
(Z_{y'},z_{y'})\in\Tan(X,x_0).
\end{equation}
Moreover, the $L$-\bilip map $f$ induces an $L$-\bilip map
\begin{equation}\label{fyprime}
f_{y'}:(Y,y')\to(Z_{y'},z_{y'}).
\end{equation}

Just a remark: as map defined on the set $Y$, $f_{y'}$ could differ from $f_{y''}$. 
The reason is that we are considering metric spaces up to isometric equivalence.
 We could make explicit the fact that $(Y,y')$ should be identified via an isometry with another tangent. 
However, such a rigor will only add heaviness on the reading.

Using the notations of \eqref{zyprime} and \eqref{fyprime},
we consider the set
$$G:=\left\{ g= f_{y''}^{-1}\circ\psi\circ f_{y'}\;:\; y',y''\in Y, \psi:(Z_{y'},z_{y'})\to (Z_{y'},z_{y''}), \psi\in\F \right\}.$$
It is immediate that $G$ is a group of $KL^2$-\bilip maps which acts transitively on $Y$.
 By taking the supremum over the $G$-orbit
of the distance function, and then the associated path metric, one gets an  $KL^2$-biLipschitz equivalent metric
with respect to which $G$ acts by isometries. Then by Montgomery-Zippin,  $G$   is a Lie group.
We conclude that any $Y$ is biLipschitz equivalent to $G/H$, where $H$ is the stabilizer of the action. 
Since the map 
$f_{y'}$ of \eqref{fyprime} is     biLipschitz, then $ Z_{y'}\in\Tan(X,x_0)$ is biLipschitz equivalent to $G/H$ as well.
 Since by assumption all tangents at $x_0$ are biLipschitz equivalent,
then they are all biLipschitz equivalent to the same $G/H$.
Finally, by biLipschitz homogeneity, all tangents at any point $x\in X$ are biLipschitz equivalent to the same $G/H$.
\qed
 \bibliography{general_bibliography}
\bibliographystyle{amsalpha}
 
 
\vskip 1in

\parbox{3.5in}{Enrico Le Donne:\\
~\\ 
Departement Mathematik\\
ETH Z\"urich\\
Ramistrasse 101\\
8092 Z\"urich\\
Switzerland\\
enrico.ledonne@math.ethz.ch}

\end{document}